\documentclass[11pt, leqno, a4paper]{amsart}
\usepackage{amsmath,amstext,amsthm,amsfonts}
\usepackage{amssymb,latexsym}
\usepackage[cp1252]{inputenc}
\usepackage[T1]{fontenc}
\usepackage{lmodern}
\usepackage{hyperref, color}
\usepackage[english]{babel}
\usepackage{mathrsfs}

\theoremstyle{plain}
\newtheorem{theorem}{\textbf{Theorem}}[section]
\newtheorem{thm}[theorem]{\textbf{Theorem}}
\newtheorem{prop}[theorem]{\textbf{Proposition}}

\def\boxit#1{\leavevmode\hbox{\vrule\vtop{\vbox{\kern.33333pt\hrule
    \kern1pt\hbox{\kern1pt\vbox{#1}\kern1pt}}\kern1pt\hrule}\vrule}}

\newcommand{\noid}{\noindent $\diamond$~}
\newcommand{\R}{\mathbb{R}}

\newcommand{\N}{\mathbb{N}}

\newcommand{\chih}{\widehat{\chi}}

\newcommand{\4}{\frac{1}{4}}
\newcommand{\2}{\frac{1}{2}}
\newcommand{\lip}{\mathrm{Lip}}

\newcommand{\Nbar}{\overline{N}}

\begin{document}

%%%%%%%%%%%%%%%%%%%%%%%%%%%%%%%%%%%%%%%%

\title[Finite index operators]{Remarks on J.~Espinar\lowercase{'s}\\``Finite index operators on
surfaces''}

\author[Pierre B\'{e}rard and Philippe Castillon]{Pierre B\'{e}rard and Philippe Castillon}

%\date{\today}
\date{\today ~[120406-berard-castillon-delta-plus-k-extension-hal.tex]}

\maketitle

%%%%%%%%%%%%%%%%%%%%%%%%%%%%%%%%%%%%%

\begin{abstract}
In this paper, we make some remarks on Jos\'{e} Espinar's paper ``Finite
index operators on surfaces'' [\texttt{arXiv:0911.3767}, to appear
in Journal of Geometric Analysis (2011)].
\end{abstract}\bigskip

\textbf{MSC}(2010): 58J50, 53A30, 53A10.\bigskip

\textbf{Keywords}: Spectral theory, positivity, minimal surface,
constant mean curvature surface.

\thispagestyle{empty}

\vspace{1.5cm}

%%%%%%%%%%%%%%%%%%%%%%%%%%%%%%%%%%%%%

\section{Introduction}\label{S-intro}

In \cite{BerCas11}, we considered operators of the form $J=\Delta +
a K - q$ on a complete non-compact Riemannian surface $(M,g)$, where
$\Delta$ is the non-negative Laplacian, and $K$ the Gaussian
curvature associated with the metric $g$. The parameter $a$ is some
positive constant, and $q$ is a non-negative locally integrable
function on $M$. More precisely, we studied the consequences, for
the geometry of the triple $(M,g;q)$, of the fact that the operator
$J$ is non-negative (in the sense of quadratic forms).
\medskip

Motivated by applications to minimal and \textsc{cmc} surfaces,
J.~Espinar \cite{Esp09} considers a different framework (see also
\cite{EsRo10}). More precisely, he considers a Riemannian surface
$(M,g)$, possibly with boundary $\partial M$ and not necessarily
complete, and operators of the form $\Delta + a K - c + P$, where
the parameters $a, c$ are positive constants, and $P$ is a
non-negative integrable function.\medskip

In this note, we consider complete surfaces without boundary, and
prove results similar to those in \cite{Esp09, EsRo10}, under weaker
assumptions. For this purpose, we apply the methods of
\cite{BerCas11}.

\section{General framework}

Generally speaking, we will use the same notations as in
\cite{BerCas11}, $(M,g)$ will denote a complete (possibly compact)
surface without boundary.\medskip

\subsection{The operators}

In this paper, we consider operators of the form,
\begin{equation}\label{E-op1}
J=\Delta + a K - q + P\,.
\end{equation}
Here $\Delta$ is the non-negative Laplacian, and $K$ the Gaussian
curvature associated with the metric $g$. We let $\mu$ denote the
Riemannian measure associated with $g$.

\noid We make the following assumptions on the operator $J$,
\begin{equation}\label{E-op2}
\left\{%
\begin{array}{ll}
a & \text{is a positive constant,}\\
q & \text{is a non-negative, locally integrable function on } M,\\
& \text{and we let~ } c = \inf_M q \ge 0,\\
P & \text{is an integrable function on~ } M,\\
& \text{and we let~ } \|P\|_1 = \int_M |P| \, d\mu.
\end{array}
\right.
\end{equation}
Note that we do not impose any sign condition on the function $P$.

\noid We say that the open geodesic ball $B(x_0,R)$ is
$J$-\emph{stable} if the operator $J$ is non-negative in the sense
of quadratic forms,
\begin{equation}\label{E-stab}
0 \le Q_J(\phi) = \int_M \big\{|d\phi|^2 + (aK-q+P)\phi^2 \big\} \,
d\mu
\end{equation}
for all $\phi$ in $\lip_0\big(B(x_0,R)\big)$, the Lipschitz
functions with compact support inside the ball.

\subsection{Volume growth assumptions}

Fix a reference point $x_0$ in $M$. We consider the following
assumptions on the volume growth on $(M,g)$.\smallskip

\noid We say that $(M,g)$ has \emph{polynomial volume growth of
order at most} $k$ if there exists a constant $C_k$ such that,
\begin{equation}\label{E-pvg}
V\big( B(x_0,R) \big) \le C_k (1+R)^k,
\end{equation}
for all $R > 0$.\smallskip

\noid We say that $(M,g)$ has $k$-\emph{subpolynomial volume growth}
if
\begin{equation}\label{E-spvg}
\limsup_{R\to \infty} \frac{V\big( B(x_0,R) \big)}{R^k} = 0.
\end{equation}

\noid We say that $(M,g)$ has \emph{subexponential volume growth} if
\begin{equation}\label{E-sevg}
\limsup_{R\to \infty} \frac{\ln \Big( V\big( B(x_0,R) \big)\Big)}{R}
= 0.
\end{equation}

For a complete surface without boundary, these definitions do not
depend on the choice of the reference point $x_0$, although the
constant $C_k$ a priori does.

\subsection{Fundamental inequalities}\label{ss-fi}

We briefly recall the notations of \cite{BerCas11}, Section~2. Given
a reference point $x_0 \in M$, we consider the open geodesic balls
$B(x_0,t)$, and their Euler-Poincar\'{e} characteristics
$\chi\big(B(x_0,t)\big)$. More precisely, we introduce the function,
$$
\chih(s) = \sup\big\{\chi\big(B(x_0,t)\big) ~|~ t \ge s\big\}.
$$
This is a non-increasing function with a sequence of
discontinuities, finite possibly empty, or infinite,
$\{t_j\}_{j=1}^{\Nbar}$, with $\Nbar \in \N \cup \{\infty\}$. Note
that this sequence depends on the choice of the reference point
$x_0$. We call $\omega_j$ the jump of the function $\chih$ at the
discontinuity $t_j$.\medskip

We call \emph{admissible} a function $\xi : [0,Q] \to \R$, which is
$C^1$ and piecewise $C^2$, with $\xi, \xi'' \ge 0$ and $\xi' \le 0$.
Let $N(Q)$ be the largest integer $n$ such that $t_n \le Q$.
\medskip

We now recall two key results from \cite{BerCas11}.\smallskip

\noid The topology of $M$ is controlled by the function $\chih$.
More precisely, we have the inequality (see \cite{BerCas11},
Lemma~2.1),
\begin{equation}\label{E-top}
1 - \sum_{n=1}^{\Nbar} \omega_n \le \chi(M).
\end{equation}

\noid Assume that the operator $J$ satisfies the assumptions
\eqref{E-op2}, and let $B(x_0,Q)$ be some $J$-stable ball in $M$.
Let $\xi$ be any admissible function on $[0,Q]$, with $\xi(Q)=0$,
and let $r$ denote the distance function to the center $x_0$ of the
ball. Plugging the function $\xi(r)$ into the quadratic form for $J$
and applying \cite{BerCas11}, Lemma~2.3, we obtain the inequality
\begin{equation}\label{E-fi1}
\left\{%
\begin{array}{ll}
\int_{B(x_0,Q)} q \xi^2(r)\, d\mu  \le & 2\pi a \xi^2(0) - 2\pi a
\sum_{j=1}^{N(Q)} \omega_n \xi^2(t_n)\\[5pt]
&  + \int_{B(x_0,Q)} P \xi^2(r) \\[5pt]
& + \int_{B(x_0,Q)} \big[ (1-2a) (\xi')^2 - 2a \xi \xi''\big](r) \,
d\mu \,,\\
\end{array}
\right.
\end{equation}
which yields the weaker inequality,
\begin{equation}\label{E-fi2}
\left\{%
\begin{array}{ll}
c \int_{B(x_0,Q)} \xi^2(r)\, d\mu  \le & 2\pi a \xi^2(0)
+ \|\xi\|_{\infty}^2 \, \|P\|_1 \\[5pt]
& + \int_{B(x_0,Q)} \big[ (1-2a) (\xi')^2 - 2a \xi \xi''\big](r) \,
d\mu \,.\\
\end{array}
\right.
\end{equation}

\section{Statements}\label{S-cm-s}

Inequality (\ref{E-fi1}) shows that the case in which the operator
$J = \Delta + a K - q + P$ is non-negative --under the assumptions
(\ref{E-op2})-- is similar to the case in which the operator $\Delta
+ a K - q$ has finite index, as treated in \cite{BerCas11},
Theorem~4.1. More precisely, we have the following result.

\begin{thm}\label{T-0}
Let $(M,g)$ be a complete Riemannian surface without boundary, and
let $J$ be the ope\-rator,
$$
J = \Delta + a K - q + P,
$$
with $q \ge 0$ locally integrable and $P$ an  integrable function.
Assume that $J \ge 0$ on $\lip_0(M)$, and that either of the
following conditions holds,
\begin{itemize}
    \item[(i)] $a > \4$, or
    \item[(ii)] $a=\4$, and $(M,g)$ has subexponential volume
    growth, or
    \item[(iii)] $a \in (0,\4)$, and $(M,g)$ has $k_a$-subpolynomial
    volume growth, with $k_a = 2 + \frac{4a}{1-4a}$.
\end{itemize}
Then, either $M$ is closed, or $(M,g)$ is non-compact with finite
topology and at most quadratic area growth. In particular, $(M,g)$
is conformally equivalent to a closed Riemannian surface with at
most finitely many points removed. Furthermore, $q$ is integrable on
$(M,g)$, and we have,
\begin{equation}\label{E-T0}
\int_M q \, d\mu \le 2\pi a \, \chi(M) + \int_M Pd\mu .
\end{equation}
\end{thm}

\textbf{Remark}. When considering an operator of the form $J=\Delta
+ aK + W$, taking $q=W_-$ and $P=W_+$, the previous result gives the
following. If either of the conditions (i), (ii) or (iii) holds, and
if $W_+$ is integrable, then $W\in L^1(M,\mu)$, $M$ has finite
conformal type, and
$$
0 \le 2\pi a \, \chi(M) + \int_M Wd\mu.
$$

The interesting case, in the present framework, is the case in which
the infimum $c$ of the function $q$ is positive. We have the
following result.

\begin{thm}\label{T-1}
Let $(M,g)$ be a complete Riemannian surface without boundary, and
let $J$ be the ope\-rator,
$$
J = \Delta + a K - q + P,
$$
with $q \ge c > 0$ locally integrable, and $P$ an integrable
function on $(M,g)$. Assume that $J \ge 0$ on $\lip_0(M)$, and that
either of the following conditions holds,
\begin{itemize}
    \item[(i)] $a > \4$, or
    \item[(ii)] $a=\4$, and $(M,g)$ has subexponential volume
    growth, or
    \item[(iii)] $a \in (0,\4)$, and $(M,g)$ has polynomial volume growth of
    degree at most $k$, for some $k$.
\end{itemize}
Then, either $M$ is closed, or $(M,g)$ is non-compact with finite
topology and finite volume. In particular, $(M,g)$ is conformally
equivalent to a closed Riemannian surface with at most finitely many
points removed. In both case, $M$ compact or non-compact,
\begin{equation}\label{E-T1}
c \, V(M,g) \le \int_M q \, d\mu \le 2\pi a \, \chi(M) + \int_M P \,
d\mu,
\end{equation}
where $V(M,g)$ is the volume of $(M,g)$.
\end{thm}\medskip

\textbf{Remark}. Under conditions (i) and (ii), this result is a
direct consequence of Theorem~\ref{T-0}. Note however that we only
need a polynomial volume growth condition in (iii), without any
bound on the degree (compare with Theorem~\ref{T-0}). This is so
because the condition that $J \ge 0$, with $c>0$, is quite strong.
One might wonder whether it is possible to weaken the growth
condition in (ii). \medskip

Theorems~\ref{T-0} and \ref{T-1} have their counterparts with the
assumption that the operator $J$ is \emph{non-negative} replaced by
the assumption that the operator $J$ has \emph{finite index}. As a
matter of fact, one can immediately reduce the former case to the
latter by using the following proposition of independent interest.

\begin{prop}\label{P-1}
Let $(M,g)$ be a complete Riemannian manifold and let $W$ be a
locally integrable function on $M$. Then the operator $\Delta + W$
has finite index if and only if there exists a locally integrable
function $P$ with compact support such that the operator $\Delta + W
+ P$ is non-negative.
\end{prop}

\section{Proofs}\label{S-cm-p}

\subsection{Proof of Theorem~\ref{T-0}}

Let us first deal with the case in which $M$ is closed. In this
case, we can use the constant function $\mathbf{1}$ in the quadratic
form associated with the operator $J$,
$$
Q_J(f) = \int_M \big( |df|^2 + (aK-q+P)f^2\big) \, d\mu
$$
and (\ref{E-T0}) follows immediately from the Gauss-Bonnet
theorem.\medskip

From now on, we assume that $(M,g)$ is complete,
non-compact.\medskip

\noid \emph{Case (i)}. Assume that $B(x_0,Q)$ is a $J$-stable ball
for some $Q$. Let $\xi(t) = (1 - t/Q)_{+}^{\alpha}$, for some
$\alpha \ge 1$. Then,
\begin{equation}\label{E-T0-0}
(1-2a) (\xi')^2 - 2a \xi \xi'' = - \frac{\alpha [(4a-1)\alpha -
2a]}{Q^2} (1 - \frac{t}{Q})^{2\alpha -2}_{+}.
\end{equation}
Choose $\alpha = \frac{2a}{4a-1}$. Apply (\ref{E-fi1}) with these
choices of $\xi$ and $\alpha$. Then,
\begin{equation}\label{E-T0-1}
\int_Mq(1 - \frac{r}{Q})^{2\alpha}_{+} + 2\pi a \sum_{n=1}^{N(Q)}
\omega_n (1 - \frac{t_n}{Q})^{2\alpha} \le 2\pi a + \int_MP(1 -
\frac{r}{Q})^{2\alpha}_{+}\,.
\end{equation}

Since $M$ is complete non-compact, and under the assumption of the
theorem, inequality (\ref{E-T0-1}) holds for all $Q > 0$, and we can
let $Q$ tend to infinity. Using the monotone convergence theorem for
the left-hand side and the dominated convergence theorem for the
right-hand side, we get
$$
\int_M q \, d\mu \le 2\pi a (1-\sum_1^{\bar{N}} \omega_n) + \int_M
Pd\mu\,,
$$
and inequality (\ref{E-T0}) follows from Lemma 2.1 in
\cite{BerCas11}. This inequality implies that the topology is finite
(with a lower bound for the Euler characteristic), and that $q$ is
integrable. To show that the surface is parabolic, we prove that the
volume growth is at most quadratic. To do so, we proceed as in
\cite{BerCas11}. From (\ref{E-fi2}) and \eqref{E-T0-0}, choosing
$\alpha$ large enough, we conclude that there exists a positive
constant $C_\alpha$ such that
$$
\frac{C_\alpha}{2^{2\alpha-2}Q^2}V\big(B(x_0,(\frac{Q}{2})\big) \le
\frac{C_\alpha}{Q^2}\int_M (1 - \frac{t}{Q})^{2\alpha -2}_{+}d\mu
\le 2\pi a + \|P\|_1\,,
$$
which concludes the proof.

\noid \emph{Case (ii)}. Assume that $B(x_0,Q)$ is a $J$-stable ball.
Take $\xi (t) = e^{-\alpha t} - e^{-\alpha Q}$ for some $\alpha >
0$. Then,
$$
(\xi')^2 - \xi \xi'' = \alpha^2 e^{-\alpha t} e^{-\alpha Q}.
$$
Applying (\ref{E-fi1}) with $a=\4$ and $\xi$ as above, gives
\begin{equation}\label{E-T0-2}
\left\{%
\begin{array}{ll}
\int_{B(x_0,Q)} q\xi^2(r) \, d\mu + \frac{\pi}{2} \sum_{n=1}^{N(Q)}
\omega_n
\xi^2(t_n) \le & \\[5pt]
\hspace{1cm} \frac{\pi}{2} \xi^2(0) + \int_M P\xi^2(r)d\mu +
\frac{\alpha^2}{2} e^{-\alpha Q} \int_{B(x_0,Q)} e^{-\alpha r}\, d\mu. &\\
\end{array}
\right.
\end{equation}

Since $M$ is complete non-compact, inequality (\ref{E-T0-2}) holds
for all $Q > 0$, and we can let $Q$ tend to infinity and argue as in
\cite{BerCas11}. The point is that the last term in the right-hand
side of (\ref{E-T0-2}) goes to zero when $Q$ tends to infinity for
any fixed $\alpha >0$, because $M$ has subexponential area growth.
Using monotone and dominated convergence theorems, it follows that
\begin{equation*}\label{E-T0-2b}
\int_M qe^{-2\alpha r} \, d\mu + \frac{\pi}{2} \sum_{n=1}^{N}
\omega_n e^{-2\alpha t_n} \le \frac{\pi}{2} + \int_M Pe^{-\alpha r}d\mu .
\end{equation*}
Letting $\alpha$ tend to zero, and using \cite{BerCas11}~Lemma 2.1,
we get inequality (\ref{E-T0}). In particular, $M$ has finite
topology and $q$ is integrable. To get quadratic area growth, we use
inequality (\ref{E-fi2}) with the test function $\xi$ given in
\cite{BerCas11} Lemma~2.4. We get the inequality
$$
\left\{%
\begin{array}{l}
\frac{1}{4R^2} \int_{B(R)} e^{2(1 - \frac{r}{2R})^2} \, d\mu \le
\frac{\pi}{2}e^2 + \|P\|_1 \\[6pt]
\hspace{2cm} + \2 \alpha^2 \beta^2 e^{-\alpha Q} \int_{C(R,Q)}
e^{-\alpha r} \, d\mu\,,
\end{array}
\right.
$$
and we let $Q$ tend to infinity to finish the proof.

\noid \emph{Case (iii)}. Assume that $B(x_0,Q)$ is a $J$-stable
ball. Take $\xi (t) = (1+\epsilon t)^{-\alpha} - (1+\epsilon
Q)^{-\alpha}$ with $\epsilon>0$ and $\alpha =\frac{2a}{1-4a}$. Then,
$$
(1-2a) (\xi')^2 - 2a \xi \xi'' =
2a\epsilon^2\alpha(\alpha +1)(1+\epsilon Q)^{-\alpha}(1+\epsilon t)^{-\alpha-2}.
$$

Applying (\ref{E-fi1}) to $\xi$ we find,
\begin{equation}\label{E-T0-3}
\left\{%
\begin{array}{ll}
& \int_{B(x_0,Q)} q \xi^2(r) \, d\mu + 2\pi a \sum_{n=1}^{N(Q)}
\omega_n \xi^2(t_n) \le \\[5pt]
&  \hspace{1cm} 2\pi a \xi^2(0) + \int_M P \xi^2(r) d\mu \\[5pt]
&  \hspace{1cm} + 2a\epsilon^2\alpha(\alpha +1) (1+\epsilon
Q)^{-\alpha} \int_{B(x_0,Q)} (1+\epsilon r)^{-\alpha-2} \, d\mu.
\end{array}%
\right.
\end{equation}

Since $M$ is complete non-compact, inequality (\ref{E-T0-3}) holds
for all $Q > 0$, we can let $Q$ tend to infinity, and argue as in
\cite{BerCas11}. The point is that the last term in the right-hand
side of (\ref{E-T0-3}) goes to zero when $Q$ tends to infinity for
any fixed $\epsilon >0$, because of the assumption on the area
growth of $M$. It follows that
\begin{equation*}\label{E-T0-3b}
\int_M q(1+\epsilon t)^{-\alpha} \, d\mu + 2 \pi a \sum_{n=1}^{N}
\omega_n (1+\epsilon t_n)^{-\alpha} \le 2 \pi a +
\int_M P (1+\epsilon t)^{-\alpha} d\mu .
\end{equation*}
Letting $\epsilon$ tend to zero and using \cite{BerCas11}~Lemma~2.1,
we get (\ref{E-T0}). In particular, $M$ has finite topology and $q$
is integrable. To get the quadratic area growth, we use inequality
(\ref{E-fi2}) and the test function $\xi$ given in \cite{BerCas11}
Lemma~2.5. We get the inequality,
$$
\left\{%
\begin{array}{l}
\frac{\alpha\beta}{R^2} \int_{B(R)} (1+\frac{r}{R})^{-2\beta-2} \, d\mu \le
2 \pi a + \|P\|_1 \\[5pt]
\hspace{1cm} + 2a\epsilon^2\alpha(\alpha +1) (1+\epsilon
Q)^{-\alpha} \int_{B(x_0,Q)} (1+\epsilon r)^{-\alpha-2} \, d\mu.
\end{array}
\right.
$$
We can conclude the proof by letting $Q$ tend to infinity. \qed
\medskip

\subsection{Proof of Theorem~\ref{T-1}}

Cases (i) and (ii) are direct consequences of Theorem \ref{T-0},
applying inequality (\ref{E-T0}) to the function $q\ge c>0$. In case
(iii), we first prove that $(M,g)$ has in fact polynomial volume
growth of degree $k$ less than $2+\frac{4a}{1-4a}$, this follows
from the assumption $c>0$.

\noid \emph{Case (iii), Preliminaries.} Assume that $B(x_0,Q)$ is a
$J$-stable ball. Take $\xi (t) = (1+\epsilon t)^{-\alpha} -
(1+\epsilon Q)^{-\alpha}$ for $\epsilon, \alpha > 0$. Then,
\begin{equation*}
\left\{%
\begin{array}{ll}
(1-2a) (\xi')^2 - 2a \xi \xi'' = \alpha \epsilon^2 [(1-4a)\alpha -
2a](1+\epsilon t)^{-2\alpha-2}& \\[5pt]
\hspace{1cm }+ 2a\epsilon^2\alpha(\alpha +1)
(1+\epsilon Q)^{-\alpha}(1+\epsilon t)^{-\alpha-2}. & \\
\end{array}
\right.
\end{equation*}

Applying (\ref{E-fi2}) to $\xi$ we find,
\begin{equation}\label{E-T1-3}
\left\{%
\begin{array}{ll}
& c \int_{B(x_0,Q)} \xi^2(r) \, d\mu \le \big( 2\pi a
+ \|P\|_1\big) \xi^2(0)\\[5pt]
& \hspace{1cm} + \epsilon^2 \alpha [(1-4a)\alpha - 2a]
\int_{B(x_0,Q)}
(1+\epsilon r)^{-2\alpha-2}\, d\mu \\[5pt]
&  \hspace{1cm} + 2a\epsilon^2\alpha(\alpha +1) (1+\epsilon
Q)^{-\alpha} \int_{B(x_0,Q)} (1+\epsilon r)^{-\alpha-2} \, d\mu.
\end{array}%
\right.
\end{equation}

Call respectively $A_2$ and $A_3$ the last two terms in the
right-hand side of the preceding inequality. \smallskip

Assume that there exists a positive constant $C_k$ such that
$V\big(B(x_0,t)\big) \le C_k (1+t)^k$, for all $t > 0$. Then,
\begin{equation}\label{E-T1-3a}
\left\{%
\begin{array}{ll}
\int_{B(x_0,Q)} (1+\epsilon r)^{-\beta} \, d\mu \le C_k (1+\epsilon
Q)^{-\beta}(1+Q)^k & \\[5pt]
\hspace{1cm} + \beta \epsilon C_k \int_0^Q (1+\epsilon
t)^{-\beta-1}(1+t)^k \, dt.& \\
\end{array}
\right.
\end{equation}
Since $(M,g)$ is complete non-compact, we can let $Q$ tend to
infinity in \eqref{E-T1-3}.

\noid \emph{Case (iii) continued.} Define $k_0$ by
\begin{equation}\label{E-T1-3b}
k_0 = \inf \{k ~\big| \exists C_k \text{~such that~}
V\big(B(x_0,t)\big) \le C_k (1+t)^k, \forall t > 0\}.
\end{equation}

\emph{Claim}: $k_0 < 2 + \frac{4a}{1-4a}$. Indeed if not, let $k_1$
be such that $k_0 < k_1 < k_0 + \2$. Choose $\alpha$ such that
$2\alpha + 2 = k_1 + \2$, and $\epsilon = 1$. Using (\ref{E-T1-3a}),
one finds that the term $A_2$ in (\ref{E-T1-3}) is uniformy bounded
when $Q$ tends to infinity. Similarly, one sees that the term $A_3$
tends to zero as $Q$ tends to infinity. It follows that for any $R >
0$, one has that
\begin{equation}\label{E-T1-3e}
c \int_{B(x_0,R)} (1+ r)^{-2\alpha} \, d\mu \le C(k_1),
\end{equation}
which implies that $c \, V\big(B(x_0,R)\big) \le C(k_1)
(1+R)^{2\alpha} \le C(k_1)(1+R)^{k_0-1}$. This contradicts the
definition of $k_0$. \medskip

Since $k_0 < 2 + \frac{4a}{1-4a}$, the assumption of
Theorem~\ref{T-0}~(iii) is satisfied and we can conclude. \qed
\medskip

\subsection{Proof of Proposition~\ref{P-1}}\hfill

\noid Assume that $\Delta + W$ has finite index on $C^1_0(M)$. Then
there exists a compact $K \subset M$ such that $\Delta + W$ is
non-negative on $C_0^1(M\setminus K)$. Take $\phi$ to be a smooth
function with compact support, such that $0 \le \phi \le 1$ and
$\phi \equiv 1$ in a compact neighborhood of $K$. Given any $\psi
\in C_0^1(M)$, write $\psi$ as $\psi = \phi \psi + (1-\phi)\psi$. An
easy computation gives,
\begin{equation}\label{E-fip-1}
\left\{%
\begin{array}{ll}
\int_M |d\psi|^2 + W\psi^2 = &\\[5pt]
\hspace{1cm} \int_M |d\big((1-\phi)\psi\big)|^2 + W\big((1-\phi)\psi\big)^2 &\\[5pt]
\hspace{1cm} + \int_M W \big( \phi^2 + 2\phi(1-\phi)\big) \psi^2 &\\[5pt]
\hspace{1cm} - \2 \int_M \psi^2 \Delta \big((1-\phi)^2\big) - \int_M
\psi^2 |d\phi|^2&\\[5pt]
\hspace{1cm} + 2 \int_M \phi (1-\2\phi)|d\psi|^2\,.&\\
\end{array}
\right.
\end{equation}
Because $\Delta + W$ is non-negative in $M\setminus K$, and because
of our choice of $\phi$, the first and fourth terms in the
right-hand side of \eqref{E-fip-1} are non-negative. The other terms
can be written as $-\int_M P \psi^2$, where the function $P$ is
defined by
\begin{equation}\label{E-fip-2}
\left\{%
\begin{array}{ll}
P := & |d\phi|^2 - \Delta \big( \phi (1 - \2 \phi)\big)\\[5pt]
& - W\phi^2 - 2 \phi (1-\phi)W.\\
\end{array}
\right.
\end{equation}

Recall that $W$ is locally integrable and that $\phi$ is smooth with
compact support. It follows that $P$ is locally integrable, with
compact support. By \eqref{E-fip-1}, the operator $\Delta + W + P$
is non-negative on $C_0^1(M)$, as stated.\medskip

\noid Assume that there exists a function $P$, which is locally
integrable with compact support, such that $\Delta + W + P$ is
non-negative on $C_0^1(M)$. Let $K$ be a compact neighborhood of the
support of $P$. Then,
$$
0 \le \int_M |d\psi|^2 + W\psi^2 + P\psi^2 = \int_M |d\psi|^2 +
W\psi^2,
$$
for any $\psi \in C_0^1(M\setminus K)$, and this means that $\Delta
+ W$ is non-negative on $C_0^1(M\setminus K)$. By a result of
B.~Devyver \cite{Dev10}, this implies that $\Delta + W$ has finite
index on $C_0^1(M)$. \qed.

%%PRIVATE

%%PRIVATE

%%%%%%%%%%%%%%%%%%%%%%%%%%%%%%%%%%%%%
%\newpage

\vspace{2cm}

\begin{small}
\begin{tabular}{l}
Pierre B\'{e}rard\\
Universit\'{e} Grenoble 1\\
Institut Fourier (\textsc{ujf-cnrs})\\
B.P. 74\\
38402 Saint Martin d'H\`{e}res Cedex\\
France\\
\verb+Pierre.Berard@ujf-grenoble.fr+\\
\end{tabular}

\bigskip

\begin{tabular}{l}
Philippe Castillon\\
Universit\'{e} Montpellier II\\
D\'{e}partement des sciences math\'{e}matiques CC 51\\
I3M (\textsc{umr 5149})\\
34095 Montpellier Cedex 5\\
France\\
\verb+Philippe.Castillon@univ-montp2.fr+\\
\end{tabular}
\end{small}

\end{document}